\documentclass[a4paper,11pt,leqno]{article}

\usepackage[latin1]{inputenc}
\usepackage[T1]{fontenc} %only using latex
\usepackage[english]{babel}
\usepackage{verbatim}
\usepackage{calligra}
\usepackage{enumerate}
\usepackage{dsfont}
\usepackage{amsfonts}
\usepackage{amsmath}
\usepackage{amsthm}
\usepackage{amssymb}
\usepackage{mathtools}
\usepackage{mathrsfs}
\usepackage{color}
\usepackage{graphicx}
\usepackage{bm}

\usepackage{tikz}

\usepackage{graphicx}		%inclusion des graphiques
\usepackage[hypcap=false]{caption}

\newcommand{\Ov}[1]{\overline{#1}}

\newcommand{\vr}{\varrho}

\newcommand{\vu}{\vc{u}}

\newcommand{\vc}[1]{{\bf #1}}

\newcommand{\Div}{{\rm div}_x}
\newcommand{\Grad}{\nabla_x}

\newcommand{\dx}{\,{\rm d} {x}}

\newcommand{\D}{{\rm d}}

\theoremstyle{definition}
\newtheorem{defi}{Definition}[section]

\theoremstyle{plane}
\newtheorem{thm}[defi]{Theorem}

\newcommand{\tsl}{\textsl}

\newcommand{\R}{\mathbb{R}}

\allowdisplaybreaks

\begin{document}

\newcommand{\fra}[1]{\textcolor{blue}{[***FF: #1 ***]}}

\title{\textsc{\Large{\textbf{Some remarks on steady solutions to \\ {the Euler system} in $\R^d$}}}}

\author{\normalsize\textsl{Francesco Fanelli}$\,^1\qquad$ and $\qquad$
\textsl{Eduard Feireisl}$\,^{2,3}$ \vspace{.5cm} \\
\footnotesize{$\,^{1}\;$ \textsc{Universit\'e de Lyon, Universit\'e Claude Bernard Lyon 1}} \\ % \vspace{.1cm} \\
{\footnotesize \it Institut Camille Jordan -- UMR 5208} \\ % \vspace{.1cm} \\
{\footnotesize 43 blvd. du 11 novembre 1918, F-69622 Villeurbanne cedex, FRANCE} \vspace{.2cm} \\
\footnotesize{$\,^{2}\;$ \textsc{Institute of Mathematics of the Academy of Sciences of the Czech Republic}} \\ % \vspace{.1cm} \\
{\footnotesize \v Zitn\' a 25, CZ-115 67 Praha 1, CZECH REPUBLIC} \vspace{.2cm} \\
\footnotesize{$\,^{3}\;$ \textsc{Institute of Mathematics, Technische Universit\"{a}t Berlin}} \\ % \vspace{.1cm} \\
{\footnotesize Stra{\ss}e des 17. Juni 136, 10623 Berlin, GERMANY} \vspace{.2cm} \\
\footnotesize{Email addresses: $\,^{1}\;$\ttfamily{fanelli@math.univ-lyon1.fr}}, $\;$
\footnotesize{$\,^{2,3}\;$\ttfamily{feireisl@math.cas.cz}}
\vspace{.2cm}
}

\date\today

\maketitle

%%%%%%%%%%%%%%%%%%%%%%%%%%%%%%%%%%%%%%%%%%%%%%%%%%%%%%%%%%%%%%%%%%%%%%%%%%%%%%%%%%%%%%%
\subsubsection*{Abstract}
{\footnotesize We show that the Euler system of gas dynamics in $\R^d$, $d=2,3$, with positive far field density and arbitrary far field entropy, admits infinitely many steady solutions with compactly 
supported velocity. The same proof yields a similar result for the incompressible Euler system with variable density. In particular, these are examples of global in time
smooth (non-trivial) solutions for the corresponding time-dependent systems.}

\paragraph*{\small 2010 Mathematics Subject Classification:}{\footnotesize 35Q31 % PDEs / Eq. of math. phys. / Euler equations
(primary);
76B03, % Fluid Mech / Incompressible inviscid fluids / Existence, uniqueness and regularity theory
76N10, %Fluid Mech / Compressible fluids and gas dynamics / Existence, uniqueness and regularity theory
35F60 % PDE / General first order / Boundary value problems for systems of non-linear first-order PDEs
(secondary).}

\paragraph*{\small Keywords: }{\footnotesize Euler system; non-homogeneous fluid; steady solution; far field conditions.}

\section{Introduction} \label{s:intro}

%\label{i}

In the recent work \cite{Gavr}, Gavrilov showed that there exists a {compactly supported} non-trivial velocity field $\vc{U}\in C^\infty_c({\R}^3;\R^3)$ solving the \emph{steady} incompressible Euler equations
\begin{equation} \label{eq:E}
\left\{\begin{array}{l}       
\vc{U}\cdot\nabla_x \vc{U}\, +\, \nabla_x P\,=\,0 \\[1ex]
{\rm div}_x\,\vc{U}\,=\,0\,.
       \end{array}
\right.
\end{equation}
A different proof of this fact, based on the Grad-Shafranov \tsl{ansatz} from plasma physics, was given by Constantin, La and Vicol \cite{CoLaVi}. In addition, the solution enjoys
a remarkable orthogonality relation:
\begin{equation} \label{eq:E2}
\vc{U} \cdot \Grad P = 0.
\end{equation}

The result of Gravilov 
{may seem surprising at the first glance as 
the Euler system \eqref{eq:E}, due to non-locality (by Biot-Savart law), 
is not expected to admit compactly supported velocity field, even when its vorticity is compactly supported.}
On the contrary, this property is known in the two-dimensional case, where a typical example is the \emph{Rankine vortex} (see Section 2.2 of \cite{Maj-Bert}):
\[
\vc{u} = (u_1, u_2)(x_1,x_2)\, ,\quad u_1 = -x_2 \Psi (|x|^2)\,,\quad u_2 = x_1 \Psi (|x|^2)\,,
\]
where $\Psi$ is an arbitrary smooth function.

Moreover, observe that, if $\vr = \vr(|x|)$ is a radially symmetric smooth function, with $\vr\geq0$, it is easy to check that 
\begin{equation} \label{inEu1}
\Div \vc{u} = 0\,,\quad \Div (\vr \vc{u}) = \Grad \vr \cdot \vu = 0 
\end{equation} 
and that, if in addition $\Psi$ vanishes in a neghbourhood of origin, it is possible to find a suitable radially symmetric function $\pi=\pi(|x|)$ such that
\begin{equation} \label{inEu2}
\vr \vu \cdot \Grad \vu = - \vr (|x|) \Psi^2 (|x|^2) x = - \Grad \pi(|x|)\,.
\end{equation}
Consequenly, the couple $[\vr, \vu]$ thus constructed solves the steady non-homogeneous incompressible Euler system \eqref{inEu1}, \eqref{inEu2}. 
Note that $\vr$ can be an arbitrary non--negative $C^1$ function, vacuum is allowed and $\vr$ may be even unbounded for $|x|\to\infty$.

%%%%%%%%%%%%%%%%%%%%%%%%%%%%%%%%%%%%%%%%%%
%%%%%%%%%%%%%%%%%%%%%%%%%%%%%%%%%%%%
\subsubsection*{Acknowledgements}
%%%%%%%%%%%%%%%%%%%%%%%%%%%%%%%%%%%%%
%%%%%%%%%%%%%%%%%%%%%%%%%%%%%%%%%%
{\small
The work of F.F. has been partially supported by the LABEX MILYON (ANR-10-LABX-0070) of Universit\'e de Lyon, within the program ``Investissement d'Avenir''
(ANR-11-IDEX-0007),  and by the projects BORDS (ANR-16-CE40-0027-01) and SingFlows (ANR-18-CE40-0027), all operated by the French National Research Agency (ANR).

The work of E.F. was partially supported by the Czech Sciences Foundation (GA\v CR), Grant Agreement
18--05974S. The Institute of Mathematics of the Academy of Sciences of the Czech Republic is supported by RVO:67985840.

}

\section{Stationary solutions of the Euler system of gas dynamics} \label{s:results}

The stationary solutions to the (complete) Euler system of gas dynamics satisfy the following system of equations:
\begin{equation} \label{i1}
\left\{\begin{array}{l}     
\Div (\vr \vu) = 0, \\ [1ex]
\Div (\vr \vu \otimes \vu) + \Grad \pi(\vr,s) = 0 \\[1ex]
\Div \big(\vr s \vu\big)\,=\,0
       \end{array}
\right.\,,\qquad\qquad\ x \in \R^d\,, \quad d=2,3\,,
\end{equation}
where $\vr$ is the mass density, $\vu$ the velocity, $s$ the entropy, and $\pi = \pi(\vr,s)$ the pressure, see \tsl{e.g.} Chapter 13 of \cite{Ben-Se} for details.
There are various possibilities for choosing the state variables for this problem. As we are interested in smooth solutions, the specific choice plays no role. For the sake of simplicity, 
we consider the equation of state of polytropic gases: 
\begin{equation} \label{EOS}
\pi (\vr, s) = \vr^{\gamma} \exp (a s)\,, \quad \gamma > 1\,, \quad a > 0\,. 
\end{equation}
More general state equations can be handled in a similar manner. 

%\subsection{Stationary solutions} \label{ss:compr}

As the problem is posed on the whole space $\R^d$, the far field conditions must be prescribed. 
It can be seen, by a straightforward modification of the argument by Chae \cite[Theorem 1.1]{Chaed}, that all solutions of \eqref{i1} with compactly supported velocity field and integrable pressure, 
\[
\pi(\vr, s) \geq 0\,, \quad \int_{\R^d} \pi(\vr, s) \ \dx < \infty\,, 
\]
must be trivial, meaning $\vu = 0$. Indeed multiplying the momentum equation in \eqref{i1} on $x$ and integrating by parts yields 
\[
\int_{\R^d} \left( \vr|\vu|^2 + \pi(\vr, s) \right) \dx = 0\,.
\] 
Accordingly, we focus on positive far field values of the density and arbitrary constant values of the entropy:
\begin{equation} \label{i2}
\vu \to 0\,, \qquad \vr \to \vr_\infty > 0\,, \qquad s \to s_\infty \in \R \qquad\qquad \mbox{ as }\qquad |x| \to \infty\,,
\end{equation}
where $\vr_\infty > 0$ and $s_\infty \in \R$ are given constants.

We remark that the same argument as the previous one can be used for the incompressible system \eqref{eq:E}, to show that the pressure $P$ satisfies
\begin{equation} \label{pres}
P(x) < P_\infty \equiv \lim_{|y| \to \infty } P(y) \ \mbox{ for any }\ x \in \mathcal{O} - \mbox{a non--empty open subset of} 
\ \R^3\,.
\end{equation}

\begin{comment}
Given $s_\infty$, we may fix the constant entropy $s = s_\infty$ converting \eqref{i1}, \eqref{i2} to the \emph{isentropic} 
Euler system 
\begin{equation} \label{i3}
\left\{\begin{array}{l}     
\Div (\vr \vu) = 0, \\ [1ex]
\Div (\vr \vu \otimes \vu) + \Grad \pi(\vr,s_\infty) = 0 
       \end{array}
\right.\,,\qquad\qquad\ x \in \R^d\,, \quad d=2,3\,,
\end{equation}
supplemented with the far field conditions
\begin{equation} \label{i4}
\vu \to 0\,, \qquad \vr \to \vr_\infty > 0 \qquad\qquad \mbox{ as }\qquad |x| \to \infty\,.
\end{equation}
We shall see below that even non--constant entropy profiles can be handled in a similar manner. 

Given a compactly supported velocity field $\vc{U}$ solving the homogeneous incompressible Euler system  \eqref{eq:E}, \eqref{eq:E2},  we may shift the pressure $P$ by an additive constant so that 
\begin{equation} \label{p1}
0< \underline{P} \leq P(x) \leq \Ov{P}\quad \mbox{ for }\quad  x \in \R^d\,,\qquad\quad P(x) = P_\infty > 0 \quad \mbox{ for }\quad |x| > R_0\,.
\end{equation}
\end{comment}

As observed by Gavrilov \cite{Gavr}, if $\vc{U}$ is a compactly supported velocity field solving the homogeneous incompressible Euler system \eqref{eq:E}, \eqref{eq:E2} with a pressure $P$,
and $\Psi \in C^1_c(\R)$, then 
$\vc{u} = \Psi (P) \vc{U}$ satisfies
\begin{equation} \label{p2}
\Div \vc{u} = 0\,, \quad \vc{u} \cdot \Grad \vc{u} + \Psi^2 (P) \Grad P = 0\,,\quad \vc{u} \cdot \Grad P = 0 \quad \mbox{in}\ \R^3\,.
\end{equation}
We focus on $\Psi \in C^1_c(\R)$ such that
\begin{equation} \label{supp}
{\rm supp} [\Psi] = (b, P_\infty)\,, \quad b \leq \inf_{x \in \R^d} P(x)\,.
\end{equation}
In view of \eqref{pres}, we have that $\vu = \Psi(P) \vc{U} \not \equiv 0$.

Now, given \eqref{p2}, we look for the density and entropy in the form $\vr = \widetilde{\vr}(P)$,
$s = \widetilde{s}(P)$,
for suitable $\widetilde{\vr},\  \widetilde{s} \in C^1(\R)$. Obviously, 
\[
\Div (\vr \vu) = \widetilde{\vr}'(P) \Grad P \cdot \vu +  \widetilde{\vr}(P) \Div \vu = 0\,, \quad
\Div (\vr s \vu) = ( \widetilde{\vr} \widetilde{s} )'(P) \Grad P \cdot \vu +  \widetilde{\vr} \widetilde{s} (P) \Div \vu = 0\,,
\]
while the momentum equation yields 
\[
\Div (\vr \vu \otimes \vu) + \widetilde{\vr} (P) \Psi^2 (P) \Grad P = 0\,.
\]
Therefore, seeing that 
\[
\Grad \pi (\vr, s) = \partial_\vr \pi\big( \widetilde{\vr }(P), \widetilde{s}(P) \big)
\widetilde{\vr}'(P) \Grad P + \partial_s \pi\big( \widetilde{\vr }(P), \widetilde{s}(P) \big) \widetilde{s}'(P) \Grad P \,,
\]
we adjust $\widetilde{\vr}$, $\widetilde{s}$ so that they solve
\begin{equation} \label{p3}
\frac{\D }{\D z} \pi\big(\widetilde{\vr} (z), \widetilde{s} (z) \big) = \widetilde{\vr} (z) \Psi^2 (z) ,\  
\widetilde{\vr} = \vr_0 \geq 0 , \ \widetilde{s} = s_0  \ \mbox{for}\ z \leq b,\ 
 \ \widetilde{\vr}(P_\infty) = \vr_\infty, \ \widetilde{s} (P_\infty) = s_\infty\,.
\end{equation}
In accordance with the pressure law \eqref{EOS}, we have $\partial_\vr \pi > 0$ and $\partial_s \pi > 0$ whenever $\vr > 0$;
whence it is easy to see that \eqref{p3} admits infinitely many different solutions. In particular, for any
given $\widetilde{\vr} \in C^3(\R)$, $\widetilde{s} \in C^3(\R)$ such that
\begin{equation} \label{solv}
\begin{split}
\widetilde{\vr} = \vr_0 > 0 \ \mbox{ for }\ z \leq b\,,\quad \widetilde{\vr} = \vr_\infty \ \mbox{ for }\ z \geq P_\infty\,,\quad
\widetilde{\vr}'(z) > 0 \ \mbox{ for }\ z \in (b, P_\infty) \\ 
\widetilde{s} = s_0 \ \mbox{ for }\ z \leq b\,,\quad \widetilde{s} = s_\infty \ \mbox{ for }\ z \geq P_\infty\,,\quad
\widetilde{s}'(z) \geq 0 \ \mbox{ for }\ z \in (b, P_\infty)\,,
\end{split}
\end{equation}
we may fix 
\[
\Psi(z) = \left( \frac{1}{\widetilde{\vr}(z)}\, \frac{\D }{\D z} \pi\big(\widetilde{\vr} (z), \widetilde{s} (z) \big) \right)^{\frac{1}{2}}\,.
\]
Thus, the triplet $\big[\vr = \widetilde{\vr}(P), s = \widetilde{s}(P), \vu = \Psi(P) \vc{U} \big]$ is the desired stationary solution. 
We have shown the following result. 

\begin{thm}[Stationary solutions for the full Euler system] \label{pT1}

Suppose that the pressure satisfies the equation of state \eqref{EOS}. Let $[\vc{U}, P]$ be a smooth solution of the system 
\eqref{eq:E}, \eqref{eq:E2} in $\R^3$ with compactly supported velocity field $\vc{U}$.

Then any triplet $\big[\vr = \widetilde{\vr}(P), s = \widetilde{s}(P), \vu = \Psi(P) \vc{U} \big]$, with $\widetilde{\vr}$, 
$\widetilde{s}$, $\Psi$ satisfying \eqref{supp}, \eqref{p3}, is a smooth solution of the full Euler system \eqref{i1}, \eqref{EOS}, with the far field 
conditions \eqref{i2}. 
In particular, for any given far field conditions $\vr_\infty > 0$, $s_\infty \in \R$, the Euler system \eqref{i1}, \eqref{EOS}, \eqref{i2} 
admits infinitely many smooth solutions satisfying 
\[
\vu = 0\,, \quad \vr = \vr_\infty\,, \quad s = s_\infty \quad \mbox{ outside a bounded ball in }\ \R^3\,.
\]

\end{thm}

The choice $\widetilde{s} = s_\infty$ in \eqref{solv} yields the result for the isentropic system. In addition, in this case
we get solutions satisfying
\[
\Div \vu = 0\,,\quad
\vu \cdot \Grad \vr = 0\,.
\]
In particular, they also solve the non-homogeneous incompressible Euler system \eqref{inEu1}, \eqref{inEu2} in $\R^3$.
Finally, we point out that a
similar construction for $d=2$ is possible, with the Rankine vortices replacing Gavrilov's solution.

To conclude, we remark that smooth stationary solutions are, of course, global in time smooth solutions of the corresponding evolutionary Euler system.

%%%%%%%%%%%%%%%%%%%%%%%%%%%%%%%%%%%%%%%%%%%%%%%%%%%%%%%%%%%%%%%%%%%%%%%%%%%%%%%%%%%%%%%%%%%%%%%%%%%%%%%%%%%%%%%%%%%%%%%%%
%%%%%%%%%%%%%%%%%%%%%%%%%%%%%%%%%%%%%%%%%%%%%%%%%%%%%%%%%%%%%%%%%%%%%%%%%%%%%%%%%%%%%%%%%%%%%%%%%%%%%%%%%%%%%%%%%%%%%%%%%
%%%%%%%%%%%%%%%%%%%%%%%%%%%%%%%%%%%%%%%%%%%%%%%%%%%%%%%%%%%%%%%%%%%%%%%%%%%%%%%%%%%%%%%%%%%%%%%%%%%%%%%%%%%%%%%%%%%%%%%%%

\def\cprime{$'$} \def\ocirc#1{\ifmmode\setbox0=\hbox{$#1$}\dimen0=\ht0
  \advance\dimen0 by1pt\rlap{\hbox to\wd0{\hss\raise\dimen0
  \hbox{\hskip.2em$\scriptscriptstyle\circ$}\hss}}#1\else {\accent"17 #1}\fi}

\end{document}